\documentclass[11pt]{article}
\usepackage[latin1]{inputenc}
\usepackage{amsmath,amssymb}
\usepackage{color}
\usepackage{xcolor}
\usepackage{multicol}
\usepackage{ulem,stmaryrd}
\usepackage{latexsym,epsfig}
\usepackage{amssymb,amsmath,amsfonts,amscd}
\usepackage{exscale}
\usepackage{ifthen}
\usepackage{longtable}
\usepackage{subfigure}
\usepackage{cite}
\usepackage{lscape}
\usepackage{soul}
\usepackage{ulem}
\usepackage{cancel}

\newtheorem{theo}{\textbf{Theorem}\ }
[section]
\newtheorem{lemma}[theo]{\textbf{Lemma}\ }

\newtheorem{prop}[theo]{\textbf{Proposition}\ }

\numberwithin{equation}{section}   

\begin{document}

	\begin{center}
		{\huge  \bf A  conditioned  local limit theorem \\}
		\vspace{3mm}
			{\huge  \bf for non-negative random matrices} 
		
		\vspace{1cm}
		
		M. Peign\'e  
		$^($\footnote{
			Institut Denis Poisson UMR 7013,  Universit\'e de Tours, Universit\'e d'Orl\'eans, CNRS  France.  email: peigne@univ-tours.fr }$^)$
		$\&$ 
		{ C. Pham} $^($\footnote{
			CERADE ESAIP, 18 rue du 8 mai 1945 - CS 80022 - 49180 St-Barth\'elemy d'Anjou ; email: dpham@esaip.org}$^), ^($\footnote{
			LAREMA UMR CNRS 6093, Angers, France. }$^)$

		\today
		
	\end{center}
	
	\centerline{\bf  \small Abstract }

	\vspace{0.2cm} 
	
	Let   $(S_n)_n$ be  the random process on $\mathbb R$ driven by the product of i.i.d. non-negative random matrices and $\tau$ its exit time from $]0, +\infty[$. 
	By using the adapted strategy initiated by D. Denisov and V. Wachtel, 
	we obtain  an asymptotic estimate and bounds of the probability that  the process $(S_k)_k$ remains non negative up to time $n$ and simultaneously  belongs to some compact set $[b, b+\ell ]\subset \mathbb R^{*+}$ at time $n$.

	\vspace{0.5cm}

	
	 AMS classification  60B15, 60F15 
	
	\section{Introduction and main results}

	\subsection{Motivation}
	
	Random walks conditioned to staying positive is a popular
	topic in probability.  In addition to their own interest, such as information about the maxima and the minima, the ladder variables  and the ladder epoch  of random walks on $\mathbb R$, they are also important in view of their applications,   for instance in queuing theory, in coding the genealogy of Galton-Watson trees	or else as models for polymers and interfaces; we refer to   \cite{CC} and references therein.
	
	The first interesting question  is to  determine  the  asymptotic behavior of the exit time  from the half line $[0, +\infty[$, and then to prove limit theorems for the process restricted to this half line or conditioned to  remain there. More precisely,  let   $(S_n)_{n\geq 1} $ be  a random walk in $\mathbb R$ whose  increments are independent with common distribution. Assume that   $(S_n)_{n\geq 1} $  is centered  and let $\tau$ be  its exit time from $[0, +\infty[$. Then,  for any $a, b,  \ell>0$, as $n \to +\infty$,
	\begin{equation}\label{locallimittheo}
	\mathbb P_{a} (\tau >n, a+ S_n \in [b, b+\ell[) \sim c \ \frac{ h^+(a) h^-(b) }{n^{3 /2}}  \ \ell,
	\end{equation}
	where  $c$ is a positive constant and  $h^+$ and $h^-$ are the renewal  functions associated to $(S_n)_{n \geq 1}$,  based on ascending and descending ladder heights (in particular these functions are positive).
	  The increments being  independent and identically distributed, the classical approach relies on the Wiener Hopf factorization and  related identities associated with the names of Baxter, Pollaczek and Spitzer; important references in the field are given by Feller and Spitzer in their   books  \cite{F},  \cite{Spitzer}.

	
	Important conceptual difficulties arise  both when the random walk $(S_n)_{n \geq 1}$ is $\mathbb R^d$-valued with $d \geq 2$ (the half line being replace by a general cone of the Euclidean space), or   when the  increments of the random walk are no longer independent.
	 As far as we know, equivalent  theory based on factorizations for these processes does not exist.  In dimension $d \geq 2$, the Wiener-Hopf factorization method works when the cone is a half space but breaks down for more general cones. Any attempt to develop a theory  of fluctuations for higher-dimensional random walks deals with the question: what would play the role of ladder epochs and ladder variables? \cite{GS}; Kingman showed in particular the impossibility  of extending Baxter  and Spitzer approach to random walks in higher dimension \cite{K}.   
	 
	 In 2055,  D. Denisov and V. Wachtel  developed a new approach to study the exit time from a cone of a random walk and several consequent limit theorems \cite{DW2015}. Their strategy,  based on the approximation  of these walks suitably normalized by a Brownian motion, with a strict control of the speed of convergence, is promising, powerful and flexible. It allows in particular to approach the random walks whose jumps are not i.i.d. 

	This flexible approach could be adapted to the quantity   $S_n(x):=  \ln \vert g_n\cdots g_1(x)\vert $,  where  $(g_k)_{k \geq 1}$ is a sequence   of i.i.d. random matrices,   $x$ is a non nul vector in $ \mathbb R^d$ and $\vert  \cdot \vert  $ is the $L_1$ norm in $\mathbb R^d$ ; this process  falls within the general framework of Markov walks on $\mathbb R$  satisfying some spectral gap assumption. The behavior of the tail of the distribution of $\tau_{x, a}:= \inf\{n \geq 1: a+\ln \vert g_n\cdots g_1(x)\vert   \le   0\}$ is   known for a few years  when the random matrices are invertible or   non-negative \cite{GLP} \cite{pham2018}.  This is extended  by I. Grama, R. Lauvergnat $\&$ E. Le Page in \cite{GLL1} to the case of  Markov walks, under a spectral gap assumption. Nevertheless,   the question of a local limit theorem for $\ln \vert g_n\cdots g_1(x)\vert $ confined in a half line still resists. In \cite{GLL2} such a statement holds  for conditioned Markov walks over a finite state space, in which case the  dual driving Markov chain also satisfies  nice spectral	gap properties; unfortunately, such a property does not hold for product of random matrices since it  is not realistic to assume that the random matrices $M_n$ act projectively on a finite set.




\subsection{Notations and assumptions}

		We endow $\mathbb R^d$  with the $L_1$ norm  $\vert  \cdot \vert  $   defined by   $\displaystyle  \vert  x\vert  := \sum_{i=1}^d \vert  x_i\vert  $ for any column vector $x=(x_i)_{1\leq i \leq d}$.

	Let $\mathcal S $ be the set of $d\times d$ matrices with positive entries. We endow  $\mathcal S $ with the standard multiplication of matrices, then the set $\mathcal S$ is a semigroup. For any $g =(g(i, j))_{1 \le i, j \le d} \in \mathcal S $,  we define $v$, endow $\vert \cdot \vert $ a norm on $\mathcal S$ and define $N$ as follows,
	\vspace{-0.1cm} 
	$$v(g) := \min_{1\leq j\leq d}\Bigl(\sum_{i=1}^d g(i, j)\Bigr);  \quad\vert  g\vert  :=  \sum_{i, j =1}^d g(i, j)  \quad  {\rm and} \quad  N(g):= \max \left({1\over v(g)}, \vert  g\vert \right) .
	$$    
	 Notice that $N(g) \geq 1$ for any $g \in \mathcal S $. 
	
%
	
	 Let $\mathcal C$ be the cone of column vectors defined by $ \mathcal C  := \{x \in \mathbb R ^d\mid \forall 1 \le i \le d, x_i \geq 0 \}$ 
	and   $\mathbb X$ be the limited   standard simplex defined by $ \mathbb X := \{ x \in \mathcal C\mid   |x| =1 \}$. 
	For any $x \in \mathcal C$, we denote by $\tilde x$ the corresponding  row vector and set 
	$\tilde{ \mathcal C} = \{\tilde x\mid x \in \mathcal C\}$ and $\tilde{ \mathbb X }= \{\tilde x\mid  x \in \mathbb X\}$.

	 We consider the following actions:
	\begin{itemize}
		\item the  linear action of $\mathcal S$ on $\mathcal C$ (resp. $\tilde {\mathcal C}$) defined by $(g, x) \mapsto gx$ 	(resp.  $(g, \tilde x) \mapsto \tilde x g)$ for any  $g \in \mathcal S $ and $x \in \mathcal C$, 
		\item the  projective action of $\mathcal S$ on $\mathbb X$  (resp. $\tilde{\mathbb X}$) defined by $(g, x) \mapsto g \cdot x := \displaystyle \frac{gx}{|gx|}$ (resp.  $(g, \tilde x) \mapsto  \tilde x \cdot g = \displaystyle\frac{\tilde  x g}{|\tilde x g|})$  for any  $g \in \mathcal S $ and $x \in \mathbb X$.
	\end{itemize}
 It is noticeable that  $0< v(g)\  \vert  x\vert  \leq \vert  gx\vert  \leq \vert  g\vert  \ \vert  x\vert$ for any $x \in \mathcal C$. 
 	
%
%
%
 For any fixed $x \in \mathbb X$ and $a \geq 0$,  we denote by $\tau_{x,a}$ the first time   the random process $(a+\ln \vert g_n\cdots g_1x\vert)_n$ becomes negative, i.e. 
 \[ \tau_{x,a}  := \min \{ n \geq 1:\, a+\ln \vert g_n\cdots g_1x\vert \leq 0\}.\]

	
	 We  impose the following  assumptions  on $\mu$.

	\vspace{2mm}	
	
	 \noindent {\bf P1} {\it  \underline{Moment assumption}:    There exists $\delta_1>0$ such that $\displaystyle \int_{\mathcal S}   \vert  \ln N(g)\vert  ^{ 2+\delta_1} \mu(dg) <+\infty$.} 
	
		\vspace{2mm}	
			
	Notice that hypothesis {\bf P1} is  weaker than the one in \cite{pham2018} where exponential moments  are  required;  the argument developed in  \cite{pham2018}   is  improved by allowing various  parameters (see \cite{PW}, {\it Proof of Theorem 1.6 (d)}).
	
	\vspace{2mm}

	 {\bf P2}  {\it  \underline{Irreducibility  assumption:} There exists no affine subspaces $A$ of $\mathbb R^d$ such that $A \cap \mathcal C$ is non-empty and bounded and invariant under the action of all elements of the support of $\mu$.}
	 
		\vspace{2mm}	
	 This assumption is  classical  in the context  of product of positive random matrices, it ensures in particular that  the central limit theorem satisfied by these products is  meaningful   since the   variance is positive  (see Corollary 3 in \cite{Hennion1997}).
	
	\vspace{2mm}	
	
	 {\bf P3} {\it There exists $B>0$ such that for  $\mu$-almost all $g$   in $\mathcal S$  and any $1\leq i, j, k, l\leq d$  
  $${g(i, j) \over B}\leq  g(k, l) \leq B \ g(i, j)  .$$ 
  }
  	
 This is a classical assumption  for product of random matrices with positive entries, it was first introduced by   H. Furstenberg and H. Kesten  \cite{FK}.


	 \vspace{2mm}	
	
	 {\bf P4} {\it  \underline{Centering} The upper Lyapunov exponent $\gamma_\mu$ is equal to $0$.}
	
 \vspace{2mm}	
	
	 {\bf P5} {\it There exists  $\delta_5>0$  such that $\mu \{g \in \mathcal S  : \forall x \in \mathbb X, \ln \vert  gx\vert   \geq \delta_5 \} >0$.}	
  \vspace{2mm}	

\noindent Condition  {\bf P5}  ensures that uniformly in   $x \in \mathcal C \setminus \{0 \} $,   the probability that  the process  $( a+\ln \vert g_n\cdots g_1x\vert)_{n \geq 1}$ remains in the half line $[0, +\infty[$ is positive.	 It is satisfied for instance when   $\mu\{g \mid v(g) >1\}>0$.

	 \vspace{2mm}

	  As it is usual in studying local probabilities, one has to distinguish between ``lattice'' and ``non lattice'' cases. The ``non lattice'' assumption   ensures that the $\mathbb R$-component of the trajectories of the Markov walk $(X_n, S_n)_{n \geq 0}$ do not live in the translation of a proper subgroup of $\mathbb  R$; in the contrary case, when $\mu$ is lattice,  a phenomenon of cyclic classes appears  which involves some complications which are not interesting in our context.  We refer to  equality  (\ref{arith}) in section  \ref{preliminaries} for a precise definition in the context of products of random matrices.

		\vspace{2mm}	
		 
	 {\bf P6} {\it \underline{Non-lattice assumption}  The measure $\mu$ is non-lattice.}

	\vspace{2mm}	
	

 The tail of the distribution of $\tau_{x, a}$ has been the subject of an extensive study in  \cite{pham2018}:  under  hypotheses {\bf P1}-{\bf P5},  there exists a positive  Borel function $V: \mathbb X\times \mathbb R^+ \to   \mathbb R^+$ such that  as $n \to +\infty$, 
		$$
		\mathbb P (\tau_{x, a}  >n) \sim \frac{2}{\sigma \sqrt{2 \pi n}}V(x, a).
		$$
		In the sequel, we also need   to consider   the process $(b-\ln \vert  \tilde x g_1\cdots g_n\vert)_n,  \tilde x \in \tilde{\mathbb X}, b \in \mathbb R^+$,  associated to  the right products $g_1\ldots g_n, n \geq 1$.  We thus  also consider the stopping time  
		\[
		\tilde \tau_{\tilde x,b}  := \min \{ n \geq 1:\, b-\ln \vert \tilde x g_1\cdots g_n\vert \leq 0\}.
		\]
		 As above,  there exists a positive  Borel function $\tilde V: \mathbb X\times \mathbb R^+ \to   \mathbb R^+$ such that  as $n \to +\infty$, 
		$$
		\mathbb P (\tilde \tau_{\tilde x, b}  >n) \sim \frac{2}{\sigma \sqrt{2 \pi n}}\tilde V(\tilde x, b).
		$$
 
		At last,   as $n \to +\infty$,  the sequence  
		$\displaystyle  \left(  {a+\ln\vert g_n\cdots g_1x\vert  \over  \sigma \sqrt n}  \right)_n$ conditioned to the event $(  \tau_{x, a}  >n)$ converges    weakly  
		towards the Rayleigh distribution on $\mathbb R^+$ whose density equals $  y\ e^{-y^2/2} 1_{\mathbb R^+}(y).$ Properties of the function $V$ are  precisely stated in   section \ref{preliminaries}.

  The  natural question   is  to get a local limit theorem for the process $(a+\ln\vert g_n\cdots g_1x\vert )_{n \geq 1}$ forced to stay positive up to time $n$, in other words to  describe the   behavior   of the quantity $\mathbb P (\tau_{x, a}>n, a+\ln\vert g_n\cdots g_1x\vert \in  [b, b+\ell])$   as $n \to +\infty$, when $a, b> 0$  and $\ell>0$.

\subsection{Main statements}

	We first state a version of the Gnedenko local limit theorem.

	\begin{theo}\label{theoGnedenkocone}
		Assume  hypotheses  {\rm {\bf P1}-{\bf P6}}.    Then,   as $n \to +\infty$,  for any  $x \in \mathbb X, a\in \mathbb R$, any $b\geq 0$ and $ \ell>0$, 
		\begin{equation*}\label{gnedenko} 
			\lim_{n \to +\infty} \Big\vert  n\mathbb P (\tau_{x,a}>n, a+\ln\vert g_n\cdots g_1x\vert\in  [b, b+\ell])-        {2 \sqrt{2\pi} \over  \sigma^2 \sqrt{ n}}  V(x,a)\   b \ e^{-b^2/2n}  \ \ell\Big\vert  =0,
		\end{equation*}
		the convergence being uniform  in $ x \in \mathbb X$ and $b \geq 0$. 
	\end{theo}
   Notice that  Theorem \ref{theoGnedenkocone} says only that this probability is $o(n^{-1})$.  
	The following theorem describes an  asymptotic behavior of  $\mathbb P(\tau_{x,a} >n,  a+\ln\vert g_n\cdots g_1x\vert \in  [b, b+\ell]).$ 
Recall  that  $\Delta =  \ln \delta $ where	$\delta$ is defined in Lemma \ref{keylem}.
	
	\begin{theo}\label{theolocal}
		Assume  hypotheses   {\rm {\bf P1}-{\bf P6}}.  There exists   positive constant   $ c, C >0$ such that, for any $x \in \mathbb X, a, b \geq 0$  and $\ell >0$,  
		\begin{equation}\label{3/2upper}
			n^{3/2} \mathbb P (\tau_{x,a}>n, a+\ln\vert g_n\cdots g_1x\vert \in  [b, b+\ell])\leq C\ V(x, a) \ \tilde V(x, b) \ \ell.
		\end{equation}
		 Furthermore,  there exists $\ell_0, \Delta >0$ such that, for $\ell> \ell_0$ and $ b \geq \Delta$, 		\begin{equation}\label{3/2lower}
				\liminf_{n \to +\infty} n^{3/2} \mathbb P (\tau_{x,a}>n, a+\ln\vert g_n\cdots g_1x\vert \in  [b, b+\ell])\geq c \ V(x, a) \ \tilde V(x, b) \ \ell.
			\end{equation}
		
	\end{theo}

	As for random walks with i.i.d. increments, it is expected that   this probability is in fact equivalent to $n^{3/2}$ up to a positive constant.  
	The argument relies on a combination of what is sometimes called ``reverse time" and ``duality" in the classical theory of random walks; roughly speaking, it  relies on the fact that, for a classical random walk  $(S_n)_{n \geq 1}$ with  i.i.d.  increments, the vectors  $(S_1, S_2, \ldots, S_n)$ and $(S_n-S_{n-1}, S_n-S_{n-2}, \ldots, S_n)$ have the same distribution.  In \cite{GLL1}, this idea has been developed in the context  of Markov walks over a Markov chain with finite state space, it  is technically much more difficult and  so far,  it escapes  from the framework of random matrix products (see the paragraph  before Lemma \ref{reverse} for  more detailed explanations). In the case of non-negative random matrices,  the difference between    $\ln \vert  g_n \ldots g_1 x\vert  $ and  $ \ln    \vert  g_n \ldots g_1\vert  $ is uniformly bounded  (see Lemma \ref{keylem}  below),  one can thus avoid the precise study of the associated dual chain\footnote{This study  would require restrictive conditions on $\mu$, for example the existence  of a density.} to obtain  the above result,  a bit  less precise but still worth of interest.

	\vspace{5mm}
	
	 {\it Notation.   Let $c$ be a strictly positive constant and $\phi, \psi$ be two functions of some variable $x$; we denote by $\phi \stackrel{c}{\preceq} \psi$ (or simply $\phi \preceq \psi$) when $\phi(x) \le c\ \psi(x)$ for any value of $x$. The notation $ \phi \stackrel{c}{\asymp} \psi$ (or simply $\phi \asymp \psi $) means $\phi \stackrel{c}{\preceq} \psi \stackrel{c}{\preceq} \phi $.}
 
	\section{Preliminaries}\label{preliminaries}

	 \subsection{The killed Markov walk  on $\mathbb X\times \mathbb R$  and its harmonic function}

We consider a sequence of i.i.d.  $\mathcal S$-valued matrices $ (g_n)_{n \geq 0}$ with common distribution  $\mu$   and denote the left and right product of matrices $ L_{n,k} := g_n \ldots g_k$ and $R_{k,n} = g_k\cdots g_n$ for any $n \geq k\geq 0 $. 
	
	We fix a $\mathbb X$-valued random variable $X_0$  and consider the Markov chain $(X_n)_{n \geq0}$  defined by $X_n^{X_0} := L_{n,1} \cdot X_0$ for any $n \geq 1$; when $X_0=x$, we set for simplicity  $X_n=X_n^x$.  Similarly,    the $\tilde{\mathbb X}$-valued   Markov chain $(\tilde X_n )_{n \geq 0}$  is defined by $\tilde X_n  :=\tilde X_0 \cdot R_{1,n}$ for any $n \geq 1$, where $\tilde X_0$ is a fixed $\tilde{\mathbb X}$-valued random variable.  
	
	Notice that  the sequence $(g_{n+1}, X_{n}^x)_{n \geq 0}$  (resp.  $(  g_{n+1}, \tilde X_n^{\tilde x})_{n \geq  0}$) is a $  \mathcal S \times  \mathbb X$ valued (resp.   $  \mathcal S \times \tilde{\mathbb X}$ valued) Markov chain with initial distribution $\mu \otimes \delta_x$ (resp.  $\mu\otimes \delta_{\tilde x}$). Their respective transition probability $P$  and $Q$ are  defined by: for any $  (g, x) \in    \mathcal S \times \mathbb{X}$ and any bounded Borel function $\varphi :    \mathcal S \times \mathbb X \to \mathbb{C}, \phi:   \mathcal S \times \tilde {\mathbb X} \to \mathbb{C} $, 
	\vspace{-0.1cm}
	\[
	P\varphi ( g, x) := \int_{\mathcal S} {\varphi ( h, g\cdot x)  \mu ({\rm d}h)}  
	\quad{\rm and} \quad 
	Q\phi ( g, \tilde x) := \int_{\mathcal S} {\varphi ( h, \tilde x\cdot g )  \mu ({\rm d}h)}.\]
	%
	%
	We denote by $(\Omega=(  \mathcal S \times \mathbb X)^{\otimes \mathbb N}, \mathcal F=\mathcal B(  \mathcal S \times \mathbb X )^{\otimes \mathbb N}, (g_{n+1}, X_n^x)_{n \geq 0}, \theta, \mathbb P_x)$  the canonical probability space associated with $(g_{n+1}, X_n^x)_{n \geq 0}$,   where $\theta$ is the classical ``shift operator'' on $(  \mathcal S \times \mathbb X)^{\otimes \mathbb N}$. Similarly 
	$(\tilde \Omega, \tilde {\mathcal F}, (g_{n+1},  \tilde {X}_n^{\tilde x})_{n \geq 0},   \tilde \theta,  {\mathbb P}_{\tilde x})$  denotes  the canonical probability space associated with $(g_{n+1}, \tilde X_n^{\tilde x})_{n \geq 0}$. 
	
	We  introduce next  the functions $\rho$ ant $\tilde \rho$ defined  for any $ g \in \mathcal S  $ and $x \in  \mathbb X$ by 
	$$
	\rho(g, x):=\ln  |gx| \quad{\rm and} \quad \tilde \rho(g, \tilde x):= \ln  |\tilde x g|.
	$$
	Notice that $gx = e^{\rho(g, x)} g \cdot x$ and that  $\rho$   satisfies the ``cocycle property":   
	\[
	\rho(gh, x)=  \rho(g, h\cdot x)+\rho(h, x), \quad \forall g, h\in \mathcal S  \mbox{ and } x \in \mathbb  X.
	\]
	This yields to  the following  basic decomposition  
	\begin{equation*} 
		\ln  \vert  L_{n,1} x \vert  =  \sum_{k=0}^{n-1}  \rho(g_{k+1}, X^x_{k}) \quad {\rm and} \quad 
		\ln  \vert \tilde x R_{1,n} \vert  =  \sum_{k=0}^{n-1}  \rho(g_{k+1}, \tilde X^{\tilde x}_{k}) .
	\end{equation*} 
	This is thus natural to introduce the following Markov walks  on $\mathbb X\times \mathbb R$ and
	$\tilde {\mathbb X}\times \mathbb R$ defined by 
	$\displaystyle S_n = S_0 + \sum_{k=0}^{n-1}\rho(g_{k+1}, X^x_{k})$  and 
	$\displaystyle \tilde S_n = \tilde S_0 - \sum_{k=0}^{n-1}  \tilde \rho(g_{k+1}, X^{\tilde x}_{k})$
	where   $S_0$ and $\tilde S_0$ are real valued random variables.
	Notice that the sequences $(X_n, S_n)_{n \geq 0}$ and $ (\tilde X_n, \tilde {\mathcal S}_n)_{n \geq 0}$ are  Markov chains  on $\mathbb X\times \mathbb R$ and $\tilde {\mathbb X}\times \mathbb R$  respectively,  with transition  probability $\widetilde P$  and $\widetilde Q$ defined by: for any $(x, a) \in \mathbb X \times \mathbb R$ 
	and any bounded Borel functions $\Phi : \mathbb X \times \mathbb R \to \mathbb{C}, \Psi : \tilde{\mathbb X} \times \mathbb R \to \mathbb{C}$, 
	\[
	{\widetilde P} \Phi ( x, a) = \int_{\mathcal S} \Phi(g \cdot x, a+\rho(g, x)) \mu ({\rm d}g)
	\quad {\rm and} \quad 
	{\widetilde Q} \Psi ( \tilde x, a) = \int_{\mathcal S} \Psi({ \tilde x}\cdot g, a-
	\tilde{\rho}(g, \tilde x))
	\mu ({\rm d}g).
	\]
	For any $(x, a)\in \mathbb{X}\times \mathbb R$, we denote by   $\mathbb{P}_{ x, a}$  the probability measure  on $(\Omega, \mathcal F)$  conditioned to the event $
	(X_{0}=x, S_0 =  a)$  and by $\mathbb{E}_{x, a}$ 
	the corresponding expectation; for simplicity,  we set $\mathbb P_{x, 0}= \mathbb P_x$ and 
	$\mathbb E_{x, 0}= \mathbb E_x$.

	Hence for any  $n \geq 1$,
	\begin{eqnarray*} \label{eqn4}
		\widetilde P^n \Phi(x, a) = \mathbb E [ \Phi(L_{n,1}\cdot x, a+\ln  \vert  L_{n,1} x\vert )]= \mathbb E_{x,a} [ \Phi( X_n,  S_n)].
	\end{eqnarray*}
 Next  we consider the restriction $\widetilde P_+$ to $\mathbb X \times \mathbb R^+$ of $\widetilde P$ defined  for any $(x,a) \in \mathbb X \times \mathbb R$ by: 
	$$
	\widetilde P _+ \Phi (x, a) = \widetilde P (\Phi {\bf 1}_{\mathbb X \times \mathbb R^+})(x,a).
	$$
	Let us emphasize that $\widetilde P_+$ may not be a Markov kernel on $\mathbb X \times \mathbb R ^+$. Furthermore, if $  \tau  := \min \{ n \geq 1:\, S_n \leq 0\}  $ is  the first time the random process $(S_n)_{n \geq 1}$ becomes non-positive, it holds  for any $(x, a) \in \mathbb X \times \mathbb{R^+}$ and any bounded Borel function $\Phi : \mathbb X \times \mathbb R \to \mathbb{C}$, 
	\begin{align*}\label{eqn5}
		\widetilde P_+\Phi(x, a)=\mathbb E_{x, a}[\Phi(X_1, S_1); \tau >1]&= \mathbb E [\Phi (g_1\cdot x, a+\ln   \vert  g_1\cdot x\vert ),   a+\ln   \vert  g_1\cdot x\vert  > 0].
	\end{align*}
A positive $\widetilde P_+$-harmonic function $V$ is any function from $\mathbb X \times \mathbb R ^+$ to $\mathbb R ^+$  satisfying $ \widetilde P_+ V =V $. We extend $V$ by setting $V(x, a) =0$ for $(x, a) \in \mathbb X \times \mathbb R^-_\ast$. In other words, the function $V$ is $\widetilde P_+$-harmonic if and only if for any $x \in \mathbb X$ and  $a\ge0$,
	\begin{eqnarray*} \label{eqn7}
		V(x, a) = \mathbb E_{x, a}  [V(X_1, S_1); \tau  >1].
	\end{eqnarray*}
	Similarly, let $  \tilde \tau  := \min \{ n \geq 1:\, \tilde S_n \leq 0\}  $ be  the first time the random process $(\tilde S_n)_{n \geq 1}$ becomes non-positive; for any $(x, b) \in \tilde {\mathbb X } \times \mathbb{R^+}$ and any bounded Borel function $\Psi : \tilde {\mathbb X} \times \mathbb R \to \mathbb{C}$, 
	\begin{eqnarray*}\label{eqn5}
		\mathbb E_{\tilde x, b}[ \Psi(\tilde X_1, \tilde S_1); \tilde \tau >1]= \mathbb E [\Psi (\tilde x\cdot g_1, b-\ln   \vert  \tilde x \cdot g_1 \vert );   b-\ln   \vert  \tilde x\cdot g_1 \vert  > 0].
	\end{eqnarray*}
	From Theorem II.1 in \cite{Hennion1997}, under conditions P1-P3 introduced below, there exists a unique probability measure $\nu$ on $\mathbb X$ such that for any bounded Borel function $\varphi$ from $\mathbb X$ to $\mathbb R$,
	$$
	(\mu \ast \nu) (\varphi)= \int_{\mathcal S} \int_{\mathbb X} \varphi(g \cdot x) \nu(dx) \mu(dg) = \int_{\mathbb X} \varphi(x) \nu(dx) = \nu (\varphi).
	$$
	Such a measure is said to be $\mu$-invariant.
	When $\mathbb E[ \vert  \ln \vert  A_1\vert  \vert ]<+\infty$,  the upper Lyapunov exponent associated with $\mu$ is finite and is expressed by
	\begin{eqnarray*} \label{eqn10}
		\gamma_\mu = \int_{\mathcal S} \int_{\mathbb X}  \rho(g, x) \nu (dx) \mu (dg) .
	\end{eqnarray*}

	We are now able to  give a precise definition of a lattice distribution $\mu$.  	We say that the measure $\mu$ is {\it lattice}  if there exists $t > 0, \epsilon \in [0, 2\pi[$ and a function $\psi : \mathbb X \to \mathbb R$   such that
	\begin{equation}\label{arith}
		\forall g \in T_\mu, \forall x \in \psi(T_\mu), \quad \exp \left\{ it\rho(g, x)-i \epsilon +i (\psi(g\cdot x)-\psi(x)) \right\}=1,
	\end{equation}
	where $T_\mu$ is the closed sub-semigroup generated by the support of $\mu$.

	It is also noticeable that the process $(X_n, S_n)_n$ is a semi-markovian  random walk on $\mathbb X\times \mathbb R$ with the strictly positive variance $\sigma^2 := \displaystyle \lim_{n \to +\infty} \frac{1}{n} \mathbb E_x [S_n^2]$, for any $x \in \mathbb X$. Condition {\bf P2} implies that $\sigma^2>0$; we refer to Theorem 5 in  \cite{Hennion1997}.

	In \cite{pham2018}, we establish the asymptotic behaviour of $\mathbb P (\tau_{x,a} >n)$ by studying the $\widetilde P_+$-harmonic function $V$. Firstly, we prove the existence of a $\widetilde P_+$-harmonic function  properties.

	\begin{prop} \label{prop1}
		Assume hypotheses  {\bf P1-P5}. Then there exists a $\widetilde P_+$-harmonic Borel function $V: \mathbb X\times \mathbb R^+ \to   \mathbb R^+$, such that  $t \mapsto V(x, t)$  is increasing on $\mathbb R^+ $ for any $x \in \mathbb X$ and satisfies the following properties: there exist $c, C >0$ and $A >0$ such that for any $x \in \mathbb  X$ and $a\geq 0$,
		\[ 
		c\vee (a-A) \leq V(x, a) \leq C\ (1+a)
		\quad {\rm and} \quad    \lim_{a \to +\infty} \frac{V(x, a)}{a} =1.
		\]
		Furthermore, 	for any $x \in \mathbb X, a \ge0$ and $n \geq 1$,
		$$
		\sqrt n \mathbb P(\tau_{x, a}  > n) \leq C\  V(x,a)
		$$
		and as $n \to +\infty$, 
		$$
		\mathbb P (\tau_{x, a}  >n) \sim \frac{2}{\sigma \sqrt{2 \pi n}}V(x, a).
		$$
		At last,   as $n \to +\infty$,  the sequence  
		$\displaystyle  \left(  {a+\ln \vert L_{n, 1}x\vert \over  \sigma \sqrt n}  \right)_n$ conditioned to $(  \tau_{x, a}  >n)$ converges    weakly  
		towards the Rayleigh distribution on $\mathbb R^+$ whose density equals $  y\ e^{-y^2/2} 1_{\mathbb R^+}(y)$,  relatively to $\mathbb P_{x,a}$ for any $x \in \mathbb X$ and $a >0$.

	\end{prop}

	\subsection{Product of  positive random matrices in $\mathcal S_\delta$}

	For any fixed $B>1$, let $\mathcal S_ B$ denote the subset of $\mathcal S$ that includes   matrices satisfying {\bf P3}.   Products of random matrices are first  studied by H. Furstenberg and H. Kesten  \cite{FK} for matrices in $\mathcal S_B$ and  then being extended   to elements of $\mathcal S$ by several authors (see \cite{Hennion1997} and references therein). The restrictive condition of H. Furstenberg and H. Kesten   considerably simplifies the study. The following statement (see \cite{FK} Lemma 2) is a key argument  in the sequel to control the asymptotic behaviour of the  norm of products of matrices of $\mathcal S_B$. Let $T_{\mathcal S_B}$ be the semi-group generated by the set $\mathcal S_B$.

 \begin{lemma}\label{keylem} 
	For any $  g \in T_{\mathcal S_B}$  and $1 \le i,j,k,l \le p,$
	\begin{eqnarray} \label{eqn21}
	g(i,j) \stackrel{B^2}{\asymp} g(k,l).
	\end{eqnarray}
		In particular, there exist $ \delta >1$ such that for any $g, h \in T_{\mathcal S_B}$ and   $ x \in \mathbb X,  \tilde y \in\tilde{\mathbb{X}}$,
		\begin{enumerate}
			\item  $ \label{eqn2.11} \vert { g x} \vert  \stackrel{\small \delta}{\asymp} \vert g \vert  \,\, \mbox{and} \,\, \vert \tilde y g\vert \stackrel{\small \delta}{\asymp} \vert g\vert $, 
				\item $\vert \widetilde ygx\vert\,\, \stackrel{\small \delta}{\asymp} \,\,\vert g\vert $,
					\item $\vert g\vert \vert h\vert  \stackrel{\small \delta}{\preceq } \vert gh\vert  \leq \vert g\vert \vert h\vert$.
		\end{enumerate}
	\end{lemma}

%
%
%
%
	As a direct consequence,  the sequence $(\ln   \vert  L_{n, 1}x\vert-\ln   \vert  L_{n, 1}\vert )_{n \geq 0}$ is bounded uniformly in $ x \in \mathbb X$.  This property   is   crucial in the sequel  in order to apply the ``reverse time" trick,  an essential argument in the proofs of our main results. 

	When studying fluctuations of  random walks  $(S_n)_{n \geq 1} $ with i.i.d. increments $Y_k$ on $\mathbb R^d, d\geq 1$, it is useful  to ``reverse  time'' as follows.  For any $1\leq k \leq n$, the   random variables   $S_n-S_k= Y_{k+1}+\ldots +Y_n$  and $S_{n-k}=Y_1+\ldots +Y_{n-k}$ have the same distribution. In the case of products of random matrices, the   cocycle property $S_n(x)= \ln \vert  L_{n, 1}(x) \vert = S_k(x)+S_{n-k}(X_k)$  is more subtle and the same argument cannot be applied directly. The fact that the $g_k$ belong to $\mathcal S_\delta$ comes to our rescue here, but the price to pay is the appearance of  the constant $\Delta=\ln \delta$ that disturbs the estimates as follows. Up  to this constant $\Delta$, we can compare the distribution of $S_n(x)-S_k(x)= :  S_{n-k} (X_k)$ to  the one of $\ln \vert  g_{k+1} \ldots g_{n}\vert $, then to the one of $\ln \vert  g_{1} \ldots g_{n-k}\vert$ and at last to the one of $\ln \vert  \tilde y g_{1} \ldots g_{n-k}\vert=: -\tilde{S}_{n-k} (y)$, for any $x, y \in \mathbb X$ (notice here that for this last quantity, the non-commutativity of the product of matrices forces us  to consider  the right  linear  action  of the matrices $R_{1, n-k}$). It is the strategy that we apply repeatedly to obtain the following result.

\noindent Recall that $\Delta = \ln \delta$ where $\delta >1$ is the constant which appears in Lemma \ref{keylem}.	\begin{lemma} \label{reverse}
		For any $x, y \in \mathbb X, a, b\geq 0$ and $ \ell>0$,
		\begin{equation}\label{reversingineq1}
			\mathbb P_{x, a}(\tau >n,  S_n  \in [b, b+\ell])\leq \mathbb P_{\tilde y, \ b+\ell+\Delta}(\tilde \tau >n, \tilde S_n \in [a, a+\ell+2\Delta]).
		\end{equation}
		Similarly,  for $a\geq \ell>2\Delta >0$  and $b \geq \Delta$,
		
		\begin{equation}\label{reversingineq2}
			\mathbb P_{x, a}(\tau >n,  S_n  \in [b, b+\ell])\geq \mathbb P_{\tilde y, \ b-\Delta}(\tilde \tau >n, \tilde S_n \in [a-\ell, a-2\Delta]).
		\end{equation}
		
		
	\end{lemma}
	Proof. We begin with the demonstration of (\ref{reversingineq1}). For  any  $n \in \mathbb N, b>0$ and $h >0$, it follows that
	\begin{align*}
		\mathbb P_{x, a}(\tau >n, \  &S_n \in [b, b+\ell])\\
		&=\mathbb P_x(a+S_1>0, \ldots, a+S_{n-1}>0,  a+S_n \in [b, b+\ell])\\
		&=\mathbb P_x(a+S_n-S_{n-1}\circ \theta>0, \ldots, a+S_n-S_1\circ \theta^{n-1}>0,  a+S_n \in [b, b+\ell])\\
		&\leq \mathbb P_x(b+\ell -S_{n-1}\circ \theta>0, \ldots, b+\ell-S_1\circ \theta^{n-1}>0,  b+\ell-S_n \in [a, a+\ell]),
	\end{align*}
	 where  $\theta$ is the shift operator and  $S_{n-k}\circ \theta^k= \ln   \vert  L_{n, k+1} X_k^x\vert  \quad \ \mathbb P_x$-a.s. for any $0\leq k \leq n-1$.  By Lemma  \ref{keylem},  the quantities $\ln   \vert  L_{n, k+1} X_k^x\vert  $ and $\ln   \vert  \tilde y L_{n, k+1}\vert $ both  belong to the interval $[\ln   \vert  L_{n, k+1}\vert  -\Delta, \ln   \vert  L_{n, k+1} \vert ]$ for any $\tilde y \in \tilde{\mathbb X}$ and $0\leq k \leq n-1$. Therefore $S_{n-k}\circ \theta^k \in [\ln   \vert  \tilde y L_{n, k+1}\vert  - \Delta ; \ln   \vert  \tilde y L_{n, k+1}\vert  + \Delta ] $ and as a result 
	\begin{align*}
		\mathbb P_{x, a}(\tau >n, \  &S_n \in [b, b+\ell])\\
		&\leq \ \mathbb P(b+\ell+\Delta- \ln   \vert  \tilde y L_{n, 2}\vert  >0, \ldots, 
		b+\ell+\Delta- \ln   \vert  \tilde y L_{n, n}\vert  >0, 
		\\
		&   \qquad   \qquad \qquad  \qquad \qquad  \qquad \qquad \qquad \qquad 
		b+\ell+\Delta- \ln   \vert  \tilde y L_{n, 1}\vert    \in [a, a+ \ell +2\Delta])
		\\
		& = \ \mathbb P(b+\ell+\Delta- \ln   \vert  \tilde y R_{1, n-1}\vert  >0, \ldots, 
		b+\ell+\Delta- \ln   \vert  \tilde y R_{1, 1}\vert  >0, 
		\\
		&   \qquad   \qquad \qquad  \qquad \qquad  \qquad \qquad \qquad \qquad 
		b+\ell+\Delta- \ln   \vert  \tilde y R_{1, n}\vert    \in [a, a+ \ell +2\Delta])
		\\
		& {\rm by \ using \ the \ fact \ that\ }  (g_1, \ldots, g_n) \ {\rm and} \ (g_n, \ldots, g_1)
		\ {\rm have \ the \ same\ distribution}
		\\
		& = \ \mathbb P_{\tilde y, \ b+\ell+\Delta}(\tilde \tau >n, \tilde S_n \in [a, a+ \ell +2\Delta]). 
	\end{align*}
	
	Similarly,  for $a> \ell>2\Delta >0$  and $b >0$, we obtain the proof of (\ref{reversingineq2}) as follows.
	\begin{align*}
		\mathbb P_{x, a}(\tau >n, \  &S_n \in [b, b+\ell])\\
		&=\mathbb P_x(a+S_1>0, \ldots, a+S_{n-1}>0,  a+S_n \in [b, b+\ell])\\
		&=\mathbb P_x(a+S_n-S_{n-1}\circ \theta>0, \ldots, a+S_n-S_1\circ \theta^{n-1}>0, b \le  a+S_n \le  b+\ell)\\
		&\geq \mathbb P_x (b -S_{n-1}\circ \theta>0, \ldots, b-S_1\circ \theta^{n-1}>0,  a-\ell \le b-S_n \le a )\\
		&\geq \mathbb P(b-\Delta- \ln   \vert  \tilde y L_{n, 2}\vert  >0, \ldots, 
		b -\Delta- \ln   \vert  \tilde y L_{n, n}\vert  >0, 
		\\
		&   \qquad   \qquad \qquad  \qquad \qquad  \qquad \qquad \qquad \qquad 
		a-\ell \le b- \Delta- \ln   \vert  \tilde y L_{n, 1}\vert    \le a -2\Delta)
		\\
		&= \mathbb P(b-\Delta- \ln   \vert  \tilde y R_{1, n-1}\vert  >0, \ldots, 
		b-\Delta- \ln   \vert  \tilde y R_{1, 1}\vert  >0, 
		\\
		&   \qquad   \qquad \qquad  \qquad \qquad  \qquad \qquad \qquad \qquad 
		b-\Delta- \ln   \vert  \tilde y R_{1, n}\vert    \in [a-\ell, a-2\Delta])
		\\
		& = \mathbb P_{\tilde y, \ b-\Delta}(\tilde{S_1} >0, \ldots, \tilde{\mathcal S}_{n-1} >0, \tilde S_n \in [a-\ell, a-2\Delta])\\ 
		& = \mathbb P_{\tilde y, \ b-\Delta}(\tilde \tau >n, \tilde S_n \in [a-\ell, a -2\Delta]). 	
	\end{align*}	
	Since $a > \ell > 2\Delta > 0$, the interval $[a-\ell, a -2\Delta]$ is not empty.


	\rightline{$\Box$}

	 \subsection{Limit  theorem for product of positive random matrices}

	In this section, we state some classical results and  preparatory lemmas, useful for  the demonstration of Theorem \ref{theoGnedenkocone} and Theorem \ref{theolocal}. The following result plays a crucial role in this article.
	\begin{theo} \label{buitheo}(\cite{Bui}, Theorem 3.2.2)
		Assume hypotheses {\bf P1-P6} hold. Then for any continuous function $u: \mathbb X\to \mathbb R$ and any continuous function with compact support $\varphi: \mathbb R \to \mathbb R$, it holds
		\[
		\lim_{n \to +\infty} \Bigl\vert   \sqrt{n} \mathbb E_{x, a}  \left[ u(X_n) \varphi (S_n) \right] - { \nu(u)\over  \sigma \sqrt{2 \pi}} \int_\mathbb R \varphi (y) e^{-(y-a)^2/2\sigma^2n}{\rm d} y
		\Bigr\vert  =0,
		\]
		where the convergence is uniform in $(x, a) \in \mathbb X\times \mathbb R$. 
	\end{theo}
	We also need other elementary estimations whose proof is detailed.
	
	\begin{lemma}  \label{matriceslocal}
		 There exist  constants $c, C>0$  such that  for every $x \in \mathbb X,  a, b, \ell >0$  and  $n \geq 1$,
		\begin{equation}\label{estimate1}
			\mathbb P_{x, a}( S_n\in[ b,  b+\ell]) \leq {c  \over \sqrt{n}}\ \ell.
		\end{equation}
		  and furthermore, for any $t>0$, if  $\vert  a-b \vert  > t \sqrt{n}$ then 
		\begin{equation}\label{estimate2}
			\mathbb P_{x, a} ( S_n\in[b, b+\ell])
			\leq {  C  \over \sqrt{n}} \ \ell\  e^{-c t^2} .
	\end{equation} \end{lemma}
	{\bf Proof.}  Assertion (\ref{estimate1}) is a consequence  of  Theorem \ref{buitheo} above. 
	Assertion (\ref{estimate2}) is a more precise overestimation than   (\ref{estimate1}) 
	for large values  of the starting point $a$, namely  when $a\succeq \sqrt{n}$,  as proved below.  
	
	We fix $h, t>0$ and let $m:= \lfloor n/2 \rfloor$ be the lower round of $n/2$. We decompose $\mathbb P_{x, a}( S_n\in[b, b+\ell])$ as follows.
	\begin{align*}
		\mathbb P_{x, a}( S_n\in[b, b+\ell])&=  \mathbb P_x(a+S_n\in[b, b+\ell])\\
		&  = \underbrace{ \mathbb P_x( a+S_n\in[b, b+\ell], \vert  S_m \vert  > t\sqrt{n}/2)}_{P_1(n, x, a, b, \ell)}   +  \underbrace{ \mathbb P_x( a+S_n\in[b, b+\ell],     \vert   S_m\vert  \leq  t\sqrt{n}/2)}_{P_2(n, x, a, b, \ell)}.
	\end{align*}
	On the one hand, from the Markov property, inequality (\ref{estimate1}) and the central limit theorem for products of random matrices \cite{Hennion1997}, there exists a strictly positive constant $c$ such that, uniformly  in  $x, a$  and $b$,
	\begin{align*}
		P_1(n, x, a, b, \ell)  & =\int_{\mathbb X\times [-t\sqrt{n}/2, t\sqrt{n}/2]^c} \mathbb P_{x' } (a+a'+S_{n-m}    \in [b, b+\ell])  \quad \mathbb P_x (X_m  \in {\rm d}x',  S_m \in {\rm d}a') 
		\\ &    \leq {c \ \ell  \over \sqrt{n-m}} \mathbb P_x (\vert  S_m\vert  > t\sqrt{n}/2)  \\
		&\preceq {e^{-c t^2}\over \sqrt{n}}.  
	\end{align*}
%
 On the other hand,  when $\vert  a-b\vert   > t \sqrt{n}$,  
	 the conditions  $\vert  S_m\vert  \le t\sqrt{n}/2 $ and $a+S_n \in [b, b+\ell]$ yield $\vert S_n-S_m\vert \geq t\sqrt{n}/2-\ell$. Hence, for fixed $h$ and $n$ large enough, 
\begin{align*}
			P_2(n, x, a, b, \ell) &\leq \mathbb P_x( a+S_n\in[b, b+\ell], \vert  S_n-S_m\vert  > t\sqrt{n}/4)\\
			&
			= \mathbb P_x( a+S_m+S_{n-m}\circ \theta^m\in[b, b+\ell],  \vert  S_{n-m}\circ \theta^m \vert   > t\sqrt{n}/4)
			\\
			&= \mathbb P_x ( a+\ln \vert  L_{m, 1} X_0 \vert + \ln   \vert  L_{n, m+1} X_m\vert \in[b, b+\ell],  \vert   L_{n, m+1}  X_m  \vert   > t\sqrt{n}/4)
			\\
			&\leq 
			\mathbb P ( a+\ln \vert  L_{m, 1} x\vert + \ln   \vert  L_{n, m+1}x'\vert \in[b-\Delta, b+\ell+\Delta],\\
			& \qquad \qquad \qquad \qquad \qquad   \vert  \ln   \vert  L_{n, m+1} x'\vert   \vert   > t\sqrt{n}/4-\Delta) 
			\ {\rm for \ any} \ x' \in \mathbb X, {\rm by  \ Lemma} \  \ref{keylem}\\
			&
			\leq 
			\int_{\{  \vert  c\vert  >t\sqrt{n}/4-\Delta\}}	 
			\underbrace{\mathbb P ( a+\ln \vert  L_{m, 1}x\vert +c\in[b-\Delta, b+\ell+\Delta])}_{\displaystyle \leq \sup_{\stackrel{y \in \mathbb X}{B\in \mathbb R} } \mathbb P_y(S_m\in [B, b+\ell+2\Delta])}
			\mathbb P(  \ln   \vert  L_{n, m+1} x'\vert     \in {\rm d}c) \\
			&\preceq {1\over \sqrt{n}} 
			\mathbb P(  \vert  \ln   \vert  L_{n, m+1} x'\vert  \vert   >t\sqrt{n}/4-\Delta)\\
			&\preceq {e^{-c t^2}\over \sqrt{n}}
	\end{align*}
uniformly in $(x, a)$ by using again (\ref{estimate1}) and the central limit theorem for product of random matrices  \cite{Hennion1997}.
	
	\rightline{$\Box$}
	
	The next statement is   analogous  to the previous lemma when the walk $(a+S_n)_n$ is forced to remain positive up to time $n$.
	\begin{lemma} \label{matriceslocalconditioned}
		There exists a    constant   $C>0$  such that for all $x \in \mathbb X, a, b\geq 0,  \ell >0$  and $n\geq 1$,
		\begin{equation}\label{estimate3} 
			\mathbb P_{x, a}(\tau  >n, S_n \in [ b, b  +\ell])
			\quad \leq\quad \   C  \  {  V(x, a)\ \ell  \over n }  .
		\end{equation} 
		Furthermore,   there exists  a constant $C>0$ such that  for  any $\ell, t>0,  n \geq 1,   a>\ell +2\Delta +  t \sqrt{n}$ and $b   > \max\{ t \sqrt{n}, \Delta \}$,   
		\begin{equation}\label{estimate4}
			\mathbb P_{x,a}(\tau   \leq n,  S_n \in[ b, b  +\ell])
			\leq C{\ell \over \sqrt{n}} e^{-c t^2}.
		\end{equation} 
		
	\end{lemma}
	{\bf Proof.} For any $1 \le m \le n$,
	\begin{align*}
		\mathbb P_{x, a}(\tau  >n,  \ S_n \in [ b, b  +\ell]) &
		\leq \mathbb P_{x, a}(\tau  >m, S_n \in [ b, b  +\ell])\\
		&= \int_{\mathbb X\times \mathbb R^*_+} \mathbb P_{x', a'}(S_{n-m} \in [b, b+\ell])\ \mathbb P_{x, a}(\tau  >m, (X_m, S_m) \in {\rm d}x' {\rm d}a ')
		\\
		& \ \preceq  \  { \mathbb P_{x, a} (\tau>m)\over \sqrt{n -m}}h  \     \quad {\rm by } \ (\ref{estimate1}) \\
		& \leq c \ {  V(x, a) \ \ell \over n}   \ {\rm by\   Proposition\ } \ref{prop1}.
	\end{align*}
	To prove assertion (\ref{estimate4}),  we work in two steps. Let $m = \lfloor n/2 \rfloor $.   
	
	\underline{Step 1}.  When $b> t \sqrt{n}  $, by using the Markov property, 
	\begin{align*}
		& \mathbb P_{x,a}(\tau   \leq m,  S_n \in[ b, b  +\ell])\\
		&\qquad =\sum_{k=1}^{m  }
		\mathbb P_{x,a}(\tau   =k,   S_n  \in[ b, b  +\ell])
		\\
		&\qquad  = \sum_{k=1}^{m  } \int_{\mathbb X\times \mathbb  R^{ -}}
		\mathbb P_{x', a'} ( S_{n-k}  \in[ b, b  +\ell])
		\quad 
		\mathbb P_{x, a}(\tau   =k,  (X_k, S_k )\in {\rm d}x' {\rm d}a')  \\
		&\qquad  \leq \max_{n-m \leq \ell \leq n} \sup_{\stackrel{x' \in \mathbb X}{\vert  a'-b\vert  >t\sqrt{n}}} \mathbb P_{x', a'}( S_\ell \in [ b, b  +\ell])  \sum_{k=1}^{m  } \int_{\mathbb X\times \mathbb  R^{ -}} P_{x, a}(\tau   =k,  (X_k, S_k )\in {\rm d}x' {\rm d}a')  
		\\ 
		&\qquad  \leq \max_{n-m \leq \ell \leq n} \sup_{\stackrel{x' \in \mathbb X}{\vert  a'-b\vert  >t\sqrt{n}}} \mathbb P_{x', a'}( S_\ell \in [ b, b  +\ell]) \  \mathbb P_{x, a} (\tau \le m) 
		\\ 
		&\qquad  \leq \max_{n-m \leq \ell \leq n} \sup_{\stackrel{x' \in \mathbb X}{\vert  a'-b\vert  >t\sqrt{n}}} \ \mathbb P_{x', a'}( S_\ell \in [ b, b  +\ell])
		\\
		& \qquad \preceq{\ell \over\sqrt{n}}\ e^{-ct^2} ,\quad {\rm for \ some \ constant }\ c>0, \ {\rm by } \ (\ref{estimate2}).
	\end{align*}
%

\underline{Step 2}.   We control the term $\mathbb P_{x, a} (m<\tau  \leq n,  S_n \in[ b, b  +\ell])$. By using the same argument to prove  (\ref{reversingineq1}), it follows that
\begin{align*}
	\mathbb P_{x, a}&(m<  \tau \leq  n, \   S_n \in [b, b+\ell]) \\
	&=\mathbb P_x(\exists k \in \{m+1, \ldots, n-1\}: a+S_k\leq 0,  a+S_n \in [b, b+\ell])\\
	&=\mathbb P (\exists k \in \{m+1, \ldots, n-1\}:  a+\ln   \vert  L_{n, 1}x\vert  - \ln   \vert  L_{n, k+1}X_k^x\vert \leq 0,  a+\ln   \vert  L_{n, 1}x\vert   \in [b, b+\ell])\\
	&\leq 
	\mathbb P (\exists k \in \{m+1, \ldots, n-1\}:  b-\Delta -\ln   \vert  \tilde y L_{n, k+1}\vert  \leq 0,  a+\ln   \vert  \tilde y L_{n, 1}\vert   \in [b-\Delta, b+\ell+\Delta])\\
	&
	=\mathbb P (\exists k \in \{m+1, \ldots, n-1\}:  b-\Delta -\ln   \vert  \tilde y R_{1, n-k}\vert  \leq 0,  a+\ln   \vert  \tilde y R_{  1, n}\vert   \in [b-\Delta, b+\ell+\Delta])\\
	& (  {\rm by \ using \ again \  the \ fact \ that\ }  (g_1, \ldots, g_n) \ {\rm and} \ (g_n, \ldots, g_1)
	\ {\rm have \ the \ same\ distribution})
	\\
	&\leq\mathbb P (\exists \ell \in \{1, \ldots, m\}:  b-\Delta -\ln   \vert  \tilde y R_{1, \ell}\vert  \leq 0,   b-\Delta -\ln   \vert  \tilde y R_{  1, n}\vert   \in [a-\ell-2\Delta, a])\\
	&= \mathbb P_{\tilde y, b-\Delta}(\tilde \tau \leq m, \tilde S_n \in [a-\ell-2\Delta, a])
	\\
	& \preceq   {\ell \over \sqrt{n}} e^{-ct^2}   ,\quad {\rm for \ some \ constant }\ c>0, 
\end{align*}
where the last inequality is obtained by applying Step 1 above to the couple $(\tilde \tau ,  \tilde S_n ) $ instead of $(\tau, S_n)$,  assured by the condition $ a-\ell -2\Delta >  t \sqrt{n}$ and  $b   >  \Delta $.

\rightline{$\Box$}

\section{Proof of Theorem \ref{theoGnedenkocone}}
We adapt the proof of Theorem 5 in \cite{DW2015} and insist on the main differences. We fix two positive constants $A$ and $\epsilon$ such that $A> 2 \epsilon >0$ and split  $\mathbb R^+$ into three intervals :
$ ]A\sqrt{n}; +\infty[,\quad ]0,  2\epsilon \sqrt n[\quad $ and $ I_{n, \epsilon, A}=[ 2\epsilon \sqrt{n}, A\sqrt{n}]$.
The proof is decomposed into three steps.

 \noindent Step 1. 
\begin{equation*}
	\lim_{A \to +\infty}\limsup_{n \to +\infty}
	\left[n   \sup_{\stackrel{x  \in  \mathbb X}{b\geq A\sqrt{n}}}\mathbb P_{x, a} ( \tau  >n, S_n \in [b, b+\ell])  \right]=0.
\end{equation*}
 Step 2.  
\begin{equation*}
	\lim_{\epsilon \to 0}\limsup_{n \to +\infty}
	\left[ n   \sup_{\stackrel{x  \in  \mathbb X}{0 < b\leq 2\epsilon\sqrt{n}}}\mathbb P_{x,a}(\tau >n, S_n \in [b, b+\ell]) \right]=0.
\end{equation*}
 Step 3. For any $A>0$,
\begin{equation*}
	\lim_{\epsilon \to 0}\limsup_{n \to +\infty} \sup_{\stackrel{x  \in  \mathbb X}{b \in I_{n, \epsilon, A}}}  \Big\vert  n\mathbb P_{x,a}(\tau >n,  S_n\in [b, b+\ell])-{2\over \sigma\sqrt {2\pi n}}\ V(x, a)\   b\  \ell\ e^{-  b ^2/2n}  \Big\vert   =0.
\end{equation*}
Theorem \ref{theoGnedenkocone} follows  by combining these three steps;   the convergence is obviously uniform over $x$.

We  set $m=\lfloor n/2 \rfloor$.  

\vspace{0.2cm}

 \noindent \underline{Proof of Step 1.} Let $a  >0$  and $b \geq A\sqrt{n}$.
We rewrite $\mathbb P_{x, a}( \tau>n,  S_n\in [b, b+\ell])$ as $P_1+P_2$, where 
\[
P_1=\mathbb P_{x,a}(\tau >n,  S_m \leq A\sqrt{m},  S_n\in [b, b+\ell])
\]
and
\[
P_2=\mathbb P_{x,a}(\tau >n,S_m >A\sqrt{m},  S_n\in [b, b+\ell]).
\]
By the Markov property, Proposition \ref{prop1} and inequality  (\ref{estimate2}), for some $c>0$,
\begin{align}\label{estimP1}
	P_1 &\leq\mathbb P_{x,a}(\tau  >m, S_m  \leq A\sqrt{m},   S_n \in  [ b, b+\ell])\notag\\
	&  	 \leq \int_{\mathbb X \times ]0, A\sqrt m] } \  \mathbb P_{x,a}(\tau  >m, X_m \in {\rm d}x', S_m \in {\rm d}a') \ \mathbb P_{x', a'}(  S_{n-m} \in [ b, b+\ell]) 
	\notag\\ 
	&  	 \leq\mathbb P_{x,a}(\tau  >m, S_m  \leq A\sqrt{m}) \ \sup_{\stackrel{x' \in \mathbb X}{0 < a' \le A\sqrt m} }\mathbb P_{x', a'}(  S_{n-m} \in [ b, b+\ell]) 
	\notag\\ 
	&   \leq\mathbb P_{x,a}(\tau  >m)\sup_{\stackrel{x' \in \mathbb X}{\vert  b-a'\vert > A\sqrt{n}/4}}\mathbb P_{x', a'}(  S_{n-m} \in [ b, b+\ell]) 
	\notag\\ 
	&   \preceq  \frac{  V(x, a) }{\sqrt n} \times  {\ell \over\sqrt{n}} \ e^{-cA^2}\quad \preceq  \quad {V(x, a)\ \over n}\ \ell\ e^{-cA^2  }.
\end{align}
Similarly, by  Proposition \ref{prop1} and (\ref{estimate1}),
\begin{align}\label{estimP2}
	P_2 &\leq \mathbb P_{x,a}(\tau >m,S_m >A\sqrt{m},  S_n\in [b, b+\ell])\notag
	\\
	&\leq  \mathbb P_{x, a} \left(   S_m  > A \sqrt{m}  \Big\vert  \tau>m \right) \ 
	\mathbb P_{x,a}(\tau  >m) \
	\sup_{(x',   a') \in \mathbb X\times \mathbb R^*_+}\mathbb P_{x', a'}\left( S_{n-m} \in [ b, b+\ell]\right)
	\notag\\
	& 	 \preceq  \mathbb P_{x, a} \left( \frac{S_m}{ \sigma \sqrt m}  > \frac{A }{ \sigma }   \Big\vert  \tau>m \right) \frac{V(x, a)}{\sqrt m} \frac{ \ell }{ \sqrt{n-m} }  \notag\\
	& \preceq \frac{V(x, a)\ \ell}{n}  \int_{A/\sigma}^{+\infty} t e^{-t^2/2}{\rm dt} 	\notag\\
	&= \frac{V(x, a)\ \ell}{n}  e^{-A^2/2\sigma^2}. \end{align}
Hence, by combining (\ref{estimP1}) and (\ref{estimP2}), it follows that
$$
\lim_{A \to +\infty} \limsup_{n \to +\infty} n\  \mathbb P_{x,a}(\tau >n,  S_n\in [b, b+\ell] )\preceq \lim_{A \to +\infty} V(x, a)\  \ell\ \left(e^{-cA^2 }+e^{-A^2/2\sigma^2}\right) = 0.
$$

\noindent  \underline{Proof of Step 2.}  Assume now  $0 < b < 2\epsilon \sqrt{n}$.  
The Markov property and Proposition  \ref{prop1}  yield
\begin{align*}
	\mathbb P _{x, a}(  \tau> & n,  \ S_n\in [b, b+\ell]) 
	\\
	& \leq \sum_{ i \in \mathbb N}
	\mathbb P_{x, a}( \tau>n,  S_m \in[i, i+1[,  S_n\in [ b, b+\ell])
	\\
	& \leq \sum_{ i \in \mathbb N}
	\mathbb P_{x, a}( \tau>m,  S_m \in[i,  i+1[) 
	\sup_{ \stackrel{x' \in \mathbb X}{a'\in[i,  i+1[}} \mathbb P_{x', a'}(\tau >n-m,  S_{n-m}  \in [ b, b+\ell]) 
	\\
	& \stackrel{ \text{by (\ref{reversingineq1})}}{\leq}   \sum_{ i \in \mathbb N}
	\mathbb P_{x, a}( \tau>m,  S_m \in[i,  i+1[) \ \mathbb P_{\tilde x, b+\ell+\Delta}( \tilde \tau  >n-m, \tilde S_{n-m} \in [i, i+\ell+2\Delta +1]) 
	\\ 
	&\stackrel{ \text{by (\ref{estimate3})}}{\leq}  C {1+a\over m}  \ \sum_{ i \in \mathbb N} \mathbb P_{\tilde x, b+\ell+\Delta}( \tilde \tau  >n-m, \tilde S_{n-m} \in [i, i+\ell+2\Delta +1]) 
	\\
	&  \preceq {V(x, a) \over n}   \mathbb P_{\tilde x, b+\ell+\Delta}( \tilde \tau  >n-m) \\
	&  \preceq {V(x, a) \over n}\times  \frac{ 1+b+\ell+\Delta}{\sqrt{n-m}}  
\quad \preceq \quad  {V(x, a) (1+2 \epsilon \sqrt{n}) \over {n^{3/2}}}.
\end{align*}
We conclude the proof of Step 2 letting $n \to+\infty$, then $\epsilon \to 0$.

 \underline{Proof of Step 3.}  We fix  $b  \in I_{n, \epsilon, A}$ and  set $  m_\epsilon=\lfloor\epsilon^3n \rfloor$. We rewrite $\mathbb P_{x, a}( \tau>n,  S_n\in [b, b+\ell])$ as  follows.
\begin{align*}
	\mathbb P_{x, a}( \tau>n, &   S_n\in [b, b+\ell])\\
	& = 
	\underbrace{
		\mathbb P_{x, a}( \tau>n,   \vert  S_{n-m_\epsilon} 
		-b\vert  >\epsilon \sqrt{n},  
		S_n\in [b, b+\ell])
	}_{\Sigma_1(n, \epsilon)}
	\\
	&\qquad \qquad \qquad \qquad+  
	\underbrace{
		\mathbb P_{x, a}( \tau>n,   \vert  S_{n-m_\epsilon} -b\vert  \leq \epsilon
		\sqrt{n},  S_n\in [b, b+\ell])
	}
	_{\Sigma_2(n, \epsilon)}
\end{align*}

For $\Sigma_1(n, \epsilon)$, by the Markov property, Proposition \ref{prop1} and (\ref{estimate2}), it follows that
\begin{align*} 
	&\Sigma_1(n, \epsilon) \notag  \\
	&\ = \int_{\mathbb X\times[b-\epsilon \sqrt{n}, b+\epsilon \sqrt{n}]^c}
	\mathbb P_{x', a'} (\tau>m_\epsilon,   S_{m_\epsilon} \in [b, b+\ell])
	\  \mathbb P_{x, a} \Bigl( \tau >n-m_\epsilon,  (X_{n-m_\epsilon}, S_{n-m_\epsilon} ) \in  {\rm d}x'   {\rm d} a' \Bigr)\  \notag  \\
	&   \leq  \sup_{\stackrel{x' \in \mathbb X}{\vert a' -b \vert > \epsilon \sqrt n }}
		\mathbb  P_{x', a'}(\tau >m_\epsilon,    S_{m_\epsilon} \in [ b, b  +\ell])
		\ \mathbb P_{x, a} ( \tau >n-m_\epsilon,  S_{n-m_\epsilon}\in[b-\epsilon \sqrt{n}, b+\epsilon \sqrt{n}]^c )     \notag 	\\
	&   \le   \sup_{\stackrel{x' \in \mathbb X}{ \vert a' -b \vert > \frac{1}{\sqrt \epsilon} \sqrt{m_\epsilon}  }} 	\mathbb  P_{x', a'}(  S_{m_\epsilon} \in [ b, b  +\ell]) \ \mathbb P_{x, a} (\tau >n-m_\epsilon)    \notag 	\\
	&  \preceq  \frac{h}{ \sqrt{m_\epsilon}} e^{-c/\epsilon}  \mathbb P_{x, a} (\tau >n-m_\epsilon)      \notag 	\\
	&  \preceq \   \frac{   V(x, a) \  \ell }{ n\epsilon^{3/2} \sqrt{1-\epsilon^3} }  e^{-c/\epsilon} , \ {\rm uniformly \ in }\ b\in I_{n, \epsilon, A}.   
\end{align*}
Therefore
\begin{eqnarray}\label{sigma_1}
	\lim_{\epsilon \to 0} \limsup_{n \to +\infty} \sup_{\stackrel{x  \in  \mathbb X}{b \in I_{n, \epsilon, A}}}	\vert n \Sigma_1 (n, \varepsilon)\vert \preceq \quad \lim_{\epsilon \to 0} \limsup_{n \to +\infty}  {   V(x, a) \  \ell   \over n\epsilon^{3/2} \sqrt{1-\epsilon^3} } e^{-c/\epsilon}=0.
\end{eqnarray}
For $\Sigma_2(n, \epsilon)$, similarly we obtain
\begin{align}\label{decompsigma2}
	&\Sigma_2(n, \epsilon) \notag \\
	&= 
	\displaystyle \int_{\mathbb X\times[b-\epsilon \sqrt{n}, b+\epsilon \sqrt{n}]}
	\mathbb P_{x', a'} (\tau>m_\epsilon,   S_{m_\epsilon} \in [b, b+\ell])
	\  \mathbb P_{x, a} \Bigl( \tau >n-m_\epsilon,  (X_{n-m_\epsilon}, S_{n-m_\epsilon} ) \in  {\rm d}x'   {\rm d} a' \Bigr)
	\notag
	\\
	& =  \Sigma'_2(n, \epsilon)  - \Sigma_2''(n, \epsilon) ,
\end{align}
where 
\[ \Sigma'_2(n, \epsilon)  :=  
\int_{\mathbb X\times[b-\epsilon \sqrt{n}, b+\epsilon \sqrt{n}]}
\mathbb P_{x', a'} (S_{m_\epsilon} \in [b, b+\ell])
\ \mathbb P_{x, a} \Bigl( \tau >n-m_\epsilon,  (X_{n-m_\epsilon}, S_{n-m_\epsilon} ) \in  {\rm d}x'   {\rm d} a' \Bigr)
\]
and 
\begin{align*}
	&\Sigma''_2(n, \epsilon) \\
	&:= 
	\int_{\mathbb X\times  [b-\epsilon \sqrt{n}, b+
		\epsilon \sqrt{n}]}
	\mathbb P_{x', a'} (\tau\leq m_\epsilon,   S_{m_\epsilon} \in [b, b+\ell])
	\ \mathbb P_{x, a} \Bigl( \tau >n-m_\epsilon,  (X_{n-m_\epsilon}, S_{n-m_\epsilon} ) \in  {\rm d}x'   {\rm d} a' \Bigr). 
\end{align*}
We first treat the second term $\Sigma''_2(n, \epsilon)$. Since $b \geq 2\epsilon \sqrt{n}$,  it holds that  $a'\geq \epsilon\sqrt{n} \ge {\sqrt{m_\epsilon\over \epsilon}}$  for any $a'\in [b-\epsilon \sqrt{n}, b+
\epsilon \sqrt{n}]$. Hence, by (\ref{estimate4}) with $t = {1 \over 2 \sqrt \epsilon }$,  for such $a'$, 
\[
\mathbb P_{x', a'}(\tau \leq m_\epsilon,    S_{m_\epsilon} \in   [b, b+\ell] )\preceq  m_\epsilon^{-1/2}  e^{-c /\epsilon }=n^{-1/2} \epsilon^{-3/2} e^{-c /\epsilon }.
\]
Consequently,
\begin{align*}
	\Sigma''_2(n, \epsilon)
	&\preceq n^{-1/2} \epsilon^{-3/2} e^{-c /\epsilon }	
	\mathbb P_{x, a}  ( \tau >n-m_\epsilon,   S_{n-m_\epsilon}
	\in [b-\epsilon \sqrt{n}, b+\epsilon \sqrt{n}])  		 \notag 	\\
	&\leq 
	n^{-1/2} \epsilon^{-3/2} e^{-c /\epsilon }  \ \mathbb P_{x, a}( \tau>n- m_\epsilon)\  \notag 	\\
	&  \preceq  n^{-1/2} \epsilon^{-3/2} e^{-c /\epsilon }  { 1+a\over  \sqrt{n - m_{\epsilon}} } \quad 	  \preceq \quad  {1+a\over n \sqrt{1- \epsilon ^3} }  \epsilon^{-3/2} e^{-c /\epsilon },
\end{align*}
which   implies that
\begin{eqnarray}\label{sigma'_2}
	\lim_{\epsilon \to 0} \limsup_{n \to +\infty} \sup_{\stackrel{x  \in  \mathbb X}{b \in I_{n, \epsilon, A}}}	\vert n \Sigma''_2 (n, \varepsilon)\vert \preceq \lim_{\epsilon\to 0}   \frac{1+a}{ \sqrt{1- \epsilon ^3} }  \epsilon^{-3/2} e^{-c /\epsilon }  \quad =0.
\end{eqnarray}
It remains to control  the term $ \Sigma'_2(n, \epsilon)$. By  Theorem \ref{buitheo},  uniformly in $b\in I_{n, \epsilon, A}$, 
\begin{align}\label{sigma'_2*}
	\Sigma'_2(n, \epsilon)&=\int_{\mathbb X\times[b-\epsilon \sqrt{n}, b+\epsilon \sqrt{n}]}
	\mathbb P_{x', a'}( S_{m_\epsilon} \in [b, b+\ell])
	\mathbb P_{x, a}( \tau >n-m_\epsilon,   (X_{n-m_\epsilon},  S_{n-m_\epsilon})   \in {\rm d}x' {\rm d} a') 	\notag 	\\ 
	&  
	=\int_{\mathbb X\times[b-\epsilon \sqrt{n}, b+\epsilon \sqrt{n}]}
	{1\over \sigma  \sqrt{2\pi  m_\epsilon}}
	e^{-(  b- a') ^2/2\sigma^2 m_\epsilon}  \ \ell\  (1 +o_n(1))  \notag \\
	&\qquad \qquad \qquad\qquad  \qquad \qquad\qquad \qquad \qquad\quad 
	\mathbb P_{x, a}( \tau >n-m_\epsilon,   (X_{n-m_\epsilon},  S_{n-m_\epsilon})   \in {\rm d}x' {\rm d} a')
	\notag\\
	&\qquad \qquad \qquad \qquad \qquad \qquad \qquad \qquad \qquad  \qquad \qquad \qquad ({\rm with} \ o_n \ {\rm uniform \ in}\ b,  a', \epsilon)
	\notag\\
	&  
	={\ell \over  \sigma   \sqrt{2\pi m_\epsilon} } \ (1 +o_n(1))
	\mathbb E_{x, a}\left[ e^{-(b-S_{n-m_\epsilon})^2/2\sigma^2m_\epsilon }; b-\epsilon \sqrt{n} \leq S_{n-m_\epsilon}\leq b+\epsilon \sqrt{n};  \tau>n-m_\epsilon \right]
	\notag\\
	& = \frac{1}{ \sigma^2 \pi} 
	 \ V(x, a)\ {\ell \over \sqrt{ m_\epsilon(n-m_\epsilon)} } \ (1 +o_n(1)) \notag\\
	&\qquad \qquad \qquad\times 
	\mathbb E_{x, a}\left[ e^{-(b-S_{n-m_\epsilon})^2/2\sigma^2m_\epsilon } 1_{ [b-\epsilon \sqrt{n},  \ b+\epsilon \sqrt{n}]}(S_{n-m_\epsilon})\slash \tau>n-m_\epsilon \right].
\end{align}
The limit theorem for $(S_n)_n$ conditioned to stay in $\mathbb R^+$ (see Proposition  \ref{prop1})  combined with the second Dini's theorem yields : for every fixed $\epsilon >0$,  as $n \to +\infty$, 
\begin{align} \label{kzcba}
	\sup_{(x, b)  \in  \mathbb X\times I_{n, \epsilon, A}} &\Big\vert 
	\mathbb E_{x, a}\left[ e^{-(b-S_{n-m_\epsilon})^2/2\sigma^2m_\epsilon } 1_{ [b-\epsilon \sqrt{n},  \ b+\epsilon \sqrt{n}]}(S_{n-m_\epsilon})\slash \tau>n-m_\epsilon\right]\notag \\
	&
	\qquad\qquad\qquad\qquad -\int_{\vert \sqrt{ 1-\epsilon^3} t-{b\over \sqrt{n}}\vert <\epsilon }t e^{-t^2/2}
	e^{-(b/\sqrt{n}-  \sqrt{ 1-\epsilon^3} t)^2/2\epsilon^3}{\rm d}t
	\Big\vert \longrightarrow 0.
\end{align} 
Since this last integral  equals $ {b\over \sqrt{n} }e^{-b^2/2n}(2\pi \epsilon)^{3/2}+o(\epsilon^{3/2})$ (see \cite{DW2015} for the details),  we obtain, combining (\ref{sigma'_2*}) and  (\ref{kzcba}),
\begin{equation}\label{sigma''_2}
	\lim_{\epsilon \to 0}\limsup_{n \to +\infty} \sup_{(x, b)  \in  \mathbb X\times I_{n, \epsilon, A}} \Big\vert 
	n \Sigma'_2(n, \epsilon)-  \frac{2 \sqrt{2 \pi}}{\sigma^2}\ V(x, a)\  {b\  \ell \over \sqrt{n}} e^{- b ^2/2n}\   
	\Big\vert =0.  
\end{equation} 
The proof of Step 3 is complete by  combining (\ref{sigma_1}), 
(\ref{decompsigma2}), (\ref{sigma'_2}) and  (\ref{sigma''_2}).

\rightline{$\Box$} 

\section{\label{sectionprrof} Proof of Theorem \ref{theolocal} }

Inequality (\ref{3/2upper}) is proved in \cite{LPP2021} Corollary 3.7.  The proof of  the lower bound  $(\ref{3/2lower})$  is based on Theorem \ref{theoGnedenkocone} and is valid for  $\ell>  2 \Delta+1$ and $b\geq \Delta$. As previously, we set $m=\lfloor n/2 \rfloor$. 
By the Markov property and (\ref{reversingineq2}),
\begin{align}\label{zeudg}
	\mathbb P _{x, a}( \tau >n,&  S_n\in  [ b, b+\ell])\notag \\
	&\geq
	\mathbb P _{x, a}( \tau >n,  S_{m} \in [\sqrt n, \sqrt{2n}],  S_n\in  [ b, b+\ell])
	\notag \\
	&\geq 
	\sum_{\stackrel{k\in \mathbb N}{\sqrt{n}\leq k\leq \sqrt{2n}-1}} 
	\mathbb P _{x, a}\Bigl( \tau >n,  k \leq S_{m} \leq k+1[, \ b \leq  S_m+S_{n-m}\circ \theta^m\leq b+\ell\Bigr)
	\notag \\
	&\geq 
	\sum_{\stackrel{k\in \mathbb N}{\sqrt{n}\leq k\leq \sqrt{2n}-1}} 
	\int_{\mathbb X\times [k, k+1]} \mathbb P _{x, a} ( \tau >m, (X_m, S_m) \in {\rm d}x'{\rm d}a')\notag \\
	&  \qquad \qquad \qquad \qquad\qquad \qquad \mathbb P_{x', a'}(\tau >n-m, b  \leq  S_{n-m} \leq b+\ell)
	\notag \\
	&
	\geq 
	\sum_{\stackrel{k\in \mathbb N}{\sqrt{n}\leq k\leq \sqrt{2n}-1}} 
	\int_{\mathbb X\times [k, k+1]} \mathbb P _{x, a}\Bigl( \tau >m, (X_m, S_m) \in {\rm d}x'{\rm d}a')\notag \\
	&  \qquad\qquad \qquad \qquad\qquad \qquad \mathbb P_{\tilde x, b-\Delta}(\tilde \tau >n-m, a'-\ell  \leq \tilde S_{n-m} \leq  a'-2 \Delta)
	\notag \\
	&
	\geq 
	\sum_{\stackrel{k\in \mathbb N}{\sqrt{n}\leq k\leq \sqrt{2n}-1}} 
	\mathbb P _{x, a}\Bigl( \tau >m, k \leq S_m\leq k+1\Bigr) \notag\\
	&  \qquad\qquad \qquad \qquad\qquad \qquad \mathbb P_{\tilde x, b-\Delta}\left(\tilde \tau >n-m, k+1-\ell\leq \tilde S_{n-m} \leq  k-2 \Delta\right)
\end{align} 
By Theorem \ref{theoGnedenkocone},  there exists a    constant  $C_0>0$ such that for any  $k \in \mathbb N$ satisfying $\sqrt{n}\leq k \leq \sqrt{2n}-1$,
\[
\liminf_{n \to +\infty}n \mathbb P _{x, a}\Bigl( \tau >m, k \leq S_m\leq k+1\Bigr)\geq C_0 
\]
and
\[
\liminf_{n \to +\infty}n  \mathbb P_{\tilde x, b-\Delta}\left(\tilde \tau >n-m, k -1\leq \tilde S_{n-m} \leq  k-2 \Delta\right)\geq C_0    (\ell-  2\Delta-1).
\]
Hence, inequality (\ref{zeudg}) yields, for $n$ large enough,  
\begin{align*}
	n^2 \mathbb P _{x, a}( \tau >n,  S_n\in  [ b, b+\ell])\geq {C_0^2\over 2} (\sqrt{2n}-\sqrt{n})   (\ell-  2\Delta-1),
\end{align*}
which   implies, for such $n$, 
\[
\mathbb P _{x, a}( \tau >n,  S_n\in  [ b, b+\ell])
\succeq { \ell-  2\Delta-1  \over n^{3/2}}.  \]
This achieves the proof of inequality (\ref{3/2lower}), taking for instance $\ell_0= 4\Delta+2$.

\rightline{$\Box$}


\end{document}